\documentclass[conference]{ieeeconf}  
\IEEEoverridecommandlockouts                              

\let\labelindent\relax
\usepackage{amsmath} 
\usepackage{amssymb}  

\usepackage{tabulary}
\usepackage{hyperref}
\usepackage[english]{babel}
\usepackage{blindtext}
\usepackage[font=footnotesize]{caption}
\usepackage[makeroom]{cancel}
\usepackage{amsfonts}
\usepackage{subfig,pgfplots}
\usepackage{amsthm}  
\usepackage{amssymb}
\usepackage{graphicx}
\usepackage{dsfont}
\usepackage{pgf,tikz}\usepackage{algorithm2e}
\RestyleAlgo{ruled}
\usepackage{algpseudocode}
\pgfplotsset{ 
	compat=newest, 
	legend style =
	{font=\footnotesize },
	label style = {font=\footnotesize},
	every tick label/.append style={font=\footnotesize}
}
\usepackage{tabularx}
\usepackage{float}
\usepackage{cite}
\usepgfplotslibrary{fillbetween}
\usetikzlibrary{shapes.multipart}
\usepackage{enumitem} 
\usepackage{bbm} 

\usepackage[strict]{changepage}

\newcommand{\tcr}{\textcolor{red}}
\newcommand{\R}{\mathbb{R}}
\newcommand{\mc}{\mathcal}
\newcommand{\ov}{\overline}
\newcommand{\eps}{\varepsilon}
\newcommand{\ben}{\begin{equation*}}
	\newcommand{\een}{\end{equation*}}
\newcommand{\be}{\begin{equation}}
	\newcommand{\ee}{\end{equation}}
\newcommand{\ba}{\begin{array}}
	\newcommand{\ea}{\end{array}}

\newtheorem{remark}{Remark}
\newtheorem{theorem}{Theorem}
\newtheorem{lemma}{Lemma}
\newtheorem{assumption}{Assumption}

\newtheorem{definition}{Definition}
\newtheorem{corollary}{Corollary}
\newtheorem{proposition}{Proposition}

\newtheorem{example}{Example}
\newenvironment{contexample}{
	\addtocounter{example}{-3} \begin{example}[continued]}{
\end{example}\addtocounter{example}{3}}

\newif\ifshowproofs
\showproofstrue 

\allowdisplaybreaks
\title{On the Stability of Dynamical Multi-Commodity Flow Networks}
\author{Davide Sipione, Giacomo Como
	\thanks{D. Sipione and G. Como are with the Department of Mathematical Sciences, Politecnico di Torino, Torino, Italy (\texttt{\{davide.sipione,giacomo.como\}@ polito.it}). G. Como is also with the Department of Automatic Control, Lund University, Sweden. This work was partially supported by the Research Project PRIN 2022 ``Extracting Essential Information and Dynamics from Complex Networks'' (Grant Agreement number 2022MBC2EZ) funded by the Italian Ministry of University and Research.}}

\begin{document}
	\maketitle
	
	\begin{abstract}
		We study a class of dynamical multi-commodity flow networks in transportation networks. These are modeled as dynamical systems describing the evolution of the densities of a number of different commodities across the cells of a transportation network. Each cell is characterized by commodity-specific increasing demand functions returning the maximum outflow of each commodity from the cell as a function of the current density of that commodity, as well as a decreasing supply function returning the total maximum inflow that is allowed in the cell as a function of the current aggregate density in the cell. Every commodity is characterized by a different routing matrix, whose entries describe the turning ratios between adjacent cells. We identify a (typically convex) capacity region: for exogenous inflow vectors  belonging to that region, we prove the existence of a locally asymptotically stable free-flow equilibrium point. Building on a contraction argument, we also provide an estimate of the basin of attraction of such free-flow equilibrium point. Finally, we analyze a simple special case showing that, when the exogenous inflow vector does not belong to the region of stability, non-free flow equilibrium points might arise. \end{abstract}
	
	\section{Introduction}
	Transportation systems have been evolving in recent years due to increasing users' heterogeneity and infrastructure reshaping to tackle environmental challenges. The traffic modeling and control problem has been an important field of study for almost over a century \cite{KynWave} and \cite{KynWave2} leading to models still widely used in recent years. This is the case of the celebrated Cell Transmission Model (CTM), \cite{DAG1} and \cite{DAG2}, that efficiently simulates traffic conditions such
	as congestion, bottlenecks and shockwaves, making it a great tool for both theoretical studies and real-world applications.
	
	In this paper, we consider a class of multi-commodity dynamical flow networks, modeled as dynamical systems driven by mass conservation laws on directed multigraphs. The nodes of the graph represent junctions and the links can be assumed as roads. Multiple commodities (representing, e.g., different types of vehicles) share the same infrastructure and interact between each other. The laws of mass conservation model how traffic volume varies on each link. Demand functions determine  the maximum total outflow for each commodity and link while supply functions are shared by commodities and model the maximum total inflow for each link. Some links  work as sources (on-ramps) whose total inflow corresponds to a constant exogenous input of traffic volume coming from the external world. To make it possible for vehicles to leave the network some link  act as sinks (off-ramps) whose total outflow simply corresponds to their demand function. 
	
	This setting is widely use in the field of transportation networks, where models tend to differ based on how congestion is handled. In fact, many allocation rules can be found in the literature: a non-FIFO rule (\cite{Stability}), FIFO rules (\cite{ConvAndRob}, \cite{Distr} and \cite{Coogan}), and mixed rules (\cite{Mixed}). The main difference between these rules appear when the supply functions are not able to accomodate the whole incoming flow and thus congestion is generated. In this paper we focus on the study of the free-flow equilibrium points. Hence, if all these model act in the same way when in free-flow, then it is possible to study them together. So, we shall take in consideration models such that in free-flow the flow between two links is simply the demand function multiplied by a \textit{turning ratio}, which is the percentage of flow directed to that particular link. 
	
	Given these considerations, it is possible to characterize the stability properties of this class of models. In fact, within the space of the exogenous inflow there exists a (convex) stability region, such that if the exogenous inflow lies strictly inside it, then the network admits a locally asymptotically stable equilibrium point. However, when this condition does not occur, two distinct situations can happen: the traffic volume of some links grows unbounded or the system still converges to an equilibrium. Although this second case is hard to characterize, we shall provide an insightful example such that this congested equilibrium point can be analytically computed.
	
	This paper provides a novel mathematical model for multicommodity dynamical flow networks which groups previously studied models (e.g., \cite{SO-DTA}) that act in the same way inside the free-flow region. The stability of this class of models is studied by generalizing results obtained in single commodity scenarios (\cite{Stability},\cite{Coogan}) and the presence of a free-flow equilibrium point is characterized. It is worth mentioning that the multi-commodity dynamical flow network model considered in this paper differs from \cite{Nilsson:2014}, where the commodities adopt dynamic routing rules on acyclic networks. Similar classes of multi-commodity flow networks are also studied  in \cite{Nilsson:2023}, within the framework of strong input-to-state stability: our approach, based on $l_1$-contraction, differs from the one in \cite{Nilsson:2023} and provides tighter results (for a more specific class of dynamical flow networks). 
	
	The rest of the paper is organized as follows: in the remainder of this section we provide some basic notation. In Section II we introduce the setting of transportation networks and we exploit it to introduce the class of models under study, providing some useful examples of different models. In Section III we study the stability of these models and state our main results. Section IV provides an insightful analytical example. Finally, in Section V conclusions are drawn and possible future research topics are issued. 
	
	\subsection*{Notation}
	The sets $\mathbb{R}$ and $\mathbb{R}_+$ represent the sets of real and non-negative numbers, respectively. Given two finite sets $\mathcal{A}$ and $\mathcal{B}$ then $\mathbb{R}^\mathcal{A}$ is the set of real numbers indexed by elements of A. Similarly, $\mathbb{R}^{\mathcal{A} \times \mathcal{B}}$ is the set indexed by the product set of $\mathcal{A}$ and $\mathcal{B}$. A directed multigraph $\mathcal{G} = (\mathcal{V},\mathcal{E})$ has a set of nodes $\mathcal{V}$ and a finite multi-set of links $\mathcal{E}$. 
	Two vectors $\sigma$ in $\mathcal{V}^\mathcal{E}$ and $\tau$ in $\mathcal{V}^\mathcal{E}$ such that $\sigma_i$ and $\tau_i$ are two nodes representing the tail and head node of a link $i$ in $ \mathcal{E}$. 
	
	\section{Model}
	
	\subsection{Multi-commodity Transportation Networks}
	We model the topology of the transportation network as a nonempty finite directed multigraph $\mathcal{G} = (\mathcal{V},\mathcal{E})$ with node set $\mc V$ and link set $\mc E$. 
	Every link $i$ in $\mc E$ is directed from its tail node $\sigma_i$ in $\mc V$ to its head node $\tau_i$ in $\mc V\setminus\{\sigma_i\}$ and  
	represents a cell (or road section). Every node $v$ in $\mc V$ represents either a junction or the interface between two adjacent cells.\footnote{Notice that assuming that $\sigma_i\ne\tau_i$ for every $i$ in $\mc E$ is equivalent to ruling out the presence of self-loops in $\mc G$.} 
	A special node $w$ in $ \mathcal{V}$ represents the external world: links $i$ in $\mc E$ with $\sigma_i = w$ are called on-ramps, while links $i$ in $\mc E$ such that $\tau_i = w$ are called off-ramps. We denote the sets of on-ramps and off-ramps as $\mathcal{R}=\{i\in\mc E:\,\sigma_i=w\}$ and $\mathcal{S}=\{i\in\mc E:\,\tau_i=w\}$, respectively. We define the set of pairs of  adjacent cells as $\mathcal{A} = \{(i,j) \in \mathcal{E} \times \mathcal{E}: \tau_i = \sigma_j \neq w\}\,.$

	\begin{figure}
		\centering
		\includegraphics[width=0.7\linewidth]{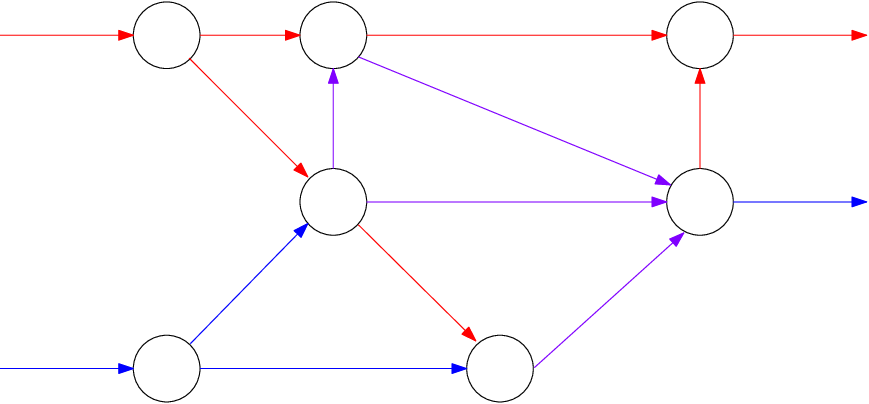}
		\caption{Two commodities share the same infrastructure: blue links are reserved for commodity 1, red links for commodity 2, and purple links can be used by both}
		\label{fig:2.1}
	\end{figure}
	\par
	
	The network infrastructure is shared by a nonempty finite set $\mathcal{K}$ of commodities. For every $i$ in $\mc E$ and $k$ in $\mc K$, the variable $x_i^{(k)}\ge0$ denotes the traffic volume (or density) of commodity $k$ on cell $i$. 
	Every non on-ramp cell $i$ in $\mathcal{E}\setminus\mc R$ is characterized by a supply function $s_i:\R_+\to\R_+$, returning the maximum possible in-flow $s_i\big(\sum_kx_i^{(k)}\big)$ as a function of the total traffic volume $\sum_kx_i^{(k)}$ on that cell, and a family of demand functions $d^{(k)}_i:\R_+\to\R_+$, for $k$ in $\mc K$, each returning the maximum out-flow $d_i^{(k)}\big(x_i^{(k)}\big)$ of commodity $k$ as a function of the traffic volume $x_i^{(k)}$ of such commodity on that cell.  Throughout the paper, we shall make the following assumption on the supply and demand functions. 
	\begin{assumption}\label{ass:demand+supply}
		\begin{enumerate}\item[(i)] The supply function $s_i$ of every non-onramp cell $i$ in $\mathcal{E}\setminus\mc R$ is Lipschitz-continuous, and such that  
			\be\label{supply-prop}s_i(0)>0\,,\qquad s_i(\xi)\le s_i(\eta)\,,\qquad \forall \xi\ge\eta\ge0\,;\ee
			\item[(ii)] the demand function $d_i^{(k)}$ of every cell $i$ in $\mc E$ and commodity $k$ in $\mc K$ is differentiable, and such that  
			\be\label{supply-prop}d_i^{(k)}(0)=0\,,\qquad (d_i^{(k)})'(\xi)>0\,,\qquad \forall \xi\ge0\,.\ee
		\end{enumerate}
	\end{assumption} 
	\begin{example}\label{ex:1}
		We introduce an examples of linear demand and affine supply function:
		\ben
		d_i^{(k)}(\xi) = \beta\xi
		\quad 
		s_i(\xi) = \gamma - \alpha\xi 
		\een
	\end{example}
	A simple example of a network with multiple commodities is shown in Figure \ref{fig:2.1}.
	
	
	To every commodity $k$ in $\mc K$, we associate a routing matrix $R^{(k)}$ in $\R_+^{\mc E\times\mc E}$, whose entries $R^{(k)}_{ij}$ represent the fraction of commodity $k$ outflow from cell $i$ that is directed to cell $j$. 
	
	\begin{assumption}\label{ass:routing-matrix}
		For every commodity $k$ in $\mc K$, the routing matrix $R^{(k)}$ in $\R_+^{\mc E\times\mc E}$ is such that 
		\be\label{R}
		R_{ij}^{(k)} = 0\,, \qquad \forall (i,j) \in \mathcal{E} \times \mathcal{E} \setminus \mathcal{A}\,,\ee
		\be\label{R2}
		\sum\limits_{j \in \mathcal{E}}{R_{ij}^{(k)}} = \left\{\ba{rl}1&\forall i \in\mc E\setminus\mathcal{S}\\0&\forall i\in\mc S \,,\ea\right.
		\ee 
		and, for every $i_0$ in $\mc E$, there exists $l\ge0$ and  $(i_1,i_2,\dots,i_l)$ in $\mc E^l$ such that 
		\begin{equation}\label{R3}
			i_l\in\mc S\,,\qquad 
			\prod_{h=1}^{l}R^{(k)}_{i_{h-1}i_h} > 0\,.
		\end{equation}
	\end{assumption} 
	Equation \eqref{R} guarantees that non-sink cells only send flow to adjacent cells, while equation \eqref{R2} implies that the outflow directed towards the external world comes only from sink cells. Moreover, equation \eqref{R3} states that for each $i$ in $\mc{E}$ there exists a finite path to a sink, so that flow can always leave the network. This assumption is rather natural and implies that, for every commodity $k$ in $\mc K$, the routing matrix $R^{(k)}$ is Schur stable, i.e., its spectral radius is smaller than $1$, so that the matrix $I-R^{(k)}$ is invertible with nonnegative inverse 
	\be\label{geometric}\big(I-R^{(k)}\big)^{-1}=I+R^{(k)}+(R^{(k)})^2+\ldots\,.\ee

	\begin{definition}
		A Multi-commodity Transportation Network (MTN) is the tuple of a nonempty finite directed multigraph $\mathcal{G} = (\mathcal{V},\mathcal{E})$, a finite set of commodities $\mathcal{K}$, supply functions $s_i$ and demand functions $d_i^{(k)}$ satisfying Assumption \ref{ass:demand+supply}, and  routing matrices $R^{(k)}$ satisfying Assumption \ref{ass:routing-matrix}, for every  cell $i$ in $\mc E$ and commodity $k$ in $\mc K$. 
	\end{definition}

	\subsection{Multi-Commodity Dynamical Flow Networks}
	
	The system's state is a time-varying element $x=x(t)$ of the nonnegative orthant 
	\be\label{X}\mc X=\R_+^{\mc E\times\mc K}\,,\ee whose entries $x_{i}^{(k)}=x_{i}^{(k)}(t)$
	represent the traffic volume  of each commodity $k$ in $\mathcal{K}$ and cell $i$ in $ \mathcal{E}$ at time $t \geq 0$. 
	We shall refer to the vector of traffic volumes of the same commodity $k$ in $\mathcal{K}$ as $x^{(k)}$ in $\R_+^{\mc E}$.  
	
	The system's inputs are  exogenous inflows $\lambda_{i}^{(k)}\geq 0$ for every on-ramp $i$ in $\mathcal{E}$ and commodity $k$ in $\mathcal{K}$,  representing the rate of traffic entering the network from the external world. For notational convenience, we shall set $\lambda_{i}^{(k)} = 0$ for every non on-ramp cell $i$ in $ \mathcal{E} \setminus \mathcal{R}$ and refer to  $\lambda=(\lambda_{i}^{(k)})_{i\in\mc E,k\in\mc K}$ as the exogenous inflow array. 
	
	Mass conservation then dictates that the variation of traffic volume of a commodity $k$  in a cell $i$ equals the difference between the total inflow to and the total outflow from that cell 
	\begin{equation}\label{eq:massCons}
		\dot{x}_{i}^{(k)}  =   \lambda_{i}^{(k)} + \sum\limits_{j \in \mathcal{E}}{f_{ji}^{(k)}(x)}  - z_{i}^{(k)}(x) \,,
	\end{equation}
	where $f_{ij}^{(k)}(x)\ge0$ is the state-dependent direct flow  of commodity $k$ in $\mathcal{K}$ from a cell $i$ in $\mc E$ to another cell $j$ in $\mc E$, while
	\begin{equation}\label{eq:zik}
		z_{i}^{(k)}(x)  =  \begin{cases}
			\sum\limits_{j \in \mathcal{E}}{f_{ij}^{(k)}(x)} & i \in \mathcal{E} \setminus \mathcal{S} \\
			d_{i}^{(k)}(x_{i}^{(k)}) & i \in \mathcal{S}\,,
		\end{cases}
	\end{equation}
	is the total out-flow of commodity $k$ from cell $i$. 
	Notice that \eqref{eq:massCons} can be rewritten in compact vector form as
	\be\label{compact}
	\dot{x}^{(k)} = \lambda^{(k)} + (R^{(k)})^{\top}z^{(k)}(x) - z^{(k)}(x) \qquad \forall k \in \mathcal{K}\,.
	\ee

	In order to fully specify the dynamics of the system, we are then left to express the functional dependence of the cell-to-cell flows $f_{ij}^{(k)}(x)$ on the state $x$. For this, we account for: 
	\begin{itemize}
		\item the demand and routing constraints, dictating that the direct flow $f_{ij}^{(k)}(x)$ of a commodity $k$ from a cell $i$ to another cell $j$ should not exceed the product of demand $d_i^{(k)}(x_i^{(k)})$ times the turning ratio $R_{ij}^{(k)}$, i.e., 
		\be\label{eq:dmnd}
		f_{ij}^{(k)}(x) \le R_{ij}^{(k)}d_{i}^{(k)}(x_{i}^{(k)}),\qquad \forall i,j\in \mathcal{E},\, k \in \mathcal{K}\,.
		\ee
		Combined with \eqref{eq:zik}, \eqref{eq:dmnd} implies that 
		$$z^{(k)}_i(x)\leq d_{i}^{(k)}(x_{i}^{(k)})\,, \qquad \forall i \in \mathcal{E}\,,\  k \in \mathcal{K}\,,   $$
		i.e., the total outflow of a commodity from  a cell never exceeds the demand, 
		while \eqref{R} and   \eqref{eq:dmnd} together imply that 
		$$f_{ij}^{(k)}(x)=0\,,\qquad \forall(i,j)\in\mc E\times\mc E\setminus\mc A\,,$$
		so that there is no direct flow of any commodity between non-adjacent cells;
		\item the supply constraints, dictating that  the aggregate total inflow of all commodities in a non on-ramp cell $i$ should not exceed the supply $s_i(\sum_kx_i^{(k)})$, i.e., 
		\be\label{eq:supp}
		\sum_{k\in\mc K}\sum\limits_{j \in \mathcal{E}}{f_{ji}^{(k)}(x)}  \leq s_i(\sum\limits_{k \in \mathcal{K}}{x_{i}^{(k)}}), \qquad i \in \mathcal{E}\setminus\mc R\,.
		\ee
		
	\end{itemize}

	It is convenient to identify those states $x$ in $\mc X$ such that  
	\be\label{free-flow} \sum\limits_{k \in \mathcal{K}}{\sum\limits_{j \in \mathcal{E}}{R_{ji}^{(k)}}d_{j}^{(k)}(x_{j}^{(k)})}  <  s_i(\sum\limits_{k \in \mathcal{K}}{x_{i}^{(k)}})\,,\ee 
	for every  $i$ in $\mathcal{E} \setminus \mathcal{R}$, i.e., 
	where the supply constraints are inactive, even if all demand and routing constraints are met with equality. We shall refer to the set of  such states as the \emph{free-flow region} and denote it by 
	$$\mathcal{F} = \{  x \in \mathcal{X}:  \eqref{free-flow}\}\,.$$ 
	
	Given all these considerations it is possible to define the Multicommodity Dynamical Flow Network.
	\begin{definition}
		Given a MTN and an exogenous inflow array $\lambda=(\lambda_{i}^{(k)})_{i\in\mc E,k\in\mc K}$, a Multicommodity Dynamical Flow Network (MDFN) is a dynamical systems satisfying equations \eqref{eq:massCons}--\eqref{eq:supp} and
		\begin{equation}\label{eq:free_const}
			f_{ij}^{(k)}(x) = R_{ij}^{(k)}d_{i}^{(k)}(x_{i}^{(k)}),\qquad \forall x\in\mc F\,,\  i,j\in \mathcal{E}\,,\  k \in \mathcal{K}\,.
		\end{equation}
	\end{definition}
	Equation $\eqref{eq:free_const}$ dictates that, for every free-flow state $x$ in $\mc F$, the demand and routing constraints \eqref{eq:dmnd} are met with equality. In particular this yields
	\begin{equation}\label{eq:free_dyn}
		z_{i}^{(k)}(x) = \sum\limits_{j \in \mathcal{E}}{f_{ij}^{(k)}(x)} = \sum\limits_{j \in \mathcal{E}}{R_{ij}^{(k)}d_{i}^{(k)}(x_{i}^{(k)})} = d_{i}^{(k)}(x_{i}^{(k)})\,,
	\end{equation}
	for every free-flow state $x$ in $\mc F$, cell $i$ in $\mathcal{E}\setminus\mc S$ and commodity $k$ in $\mathcal{K}$.
	
	Notice that the definition above only models the behavior of the system within the free-flow region, so it is possible to use models that act differently outside of that region. To this end, we shall introduce two relevant examples.
	\begin{example}
		Consider the following FIFO allocation rule
		\begin{equation}
			f_{ij}^{(k)}(x) = \gamma_i^F(x)R_{ij}^{(k)}d_{i}^{(k)}(x_{i}^{(k)})
		\end{equation}
		where $\gamma_{i}^{F}  =\max \{  \gamma \in [0,1] :    \eqref{eq:4.5}\}$, 
		\be\label{eq:4.5}
		\gamma \cdot \underset{\underset{R_{ij}^{(k)}>0}{j \in \mathcal{E}}}{\max} \:\sum\limits_{k \in \mathcal{K}}\sum\limits_{h \in \mathcal{E}} R_{hj}^{(k)}d_{h}^{(k)}(x_{h}^{(k)}) 
		\le  s_j(\sum\limits_{k \in \mathcal{K}}{x_{j}^{(k)}})\,,\ee
		i.e., $\gamma_i^F$ is the maximum value in $[0,1]$ such that for each cell $j$, with $(i,j)$ in $\mathcal{A}$, the supply constraint is satisfied. This means that all the flows sent to downstream cells are scaled by the same value, even if for some cells more flow could be allocated.
		In the literature we found examples of FIFOs models in \cite{ConvAndRob}, \cite{ConvAndRob2}, \cite{Distr}.
	\end{example}
	
	\begin{example}
		Consider the following non-FIFO allocation rule
		\begin{equation}\label{eq:10}
			f_{ij}^{(k)}(x) = R_{ij}^{(k)}d_{i}^{(k)}(x_{i}^{(k)})\min\left\{1,\frac{s_j\Big(\sum\limits_{h \in \mathcal{K}}{x_{j}^{(h)}}\Big)}{\sum\limits_{l \in \mathcal{E}}{\sum\limits_{h\in \mathcal{K}}}{R_{lj}^{(h)}d_{l}^{(h)}(x_{l}^{(h)})}} \right\}
		\end{equation}
		\normalsize
		which states that when the supply function cannot accomodate all the incoming flows of all commodities, each flow is properly scaled based on the total amount each cell and commodity want to send.
		In \cite{Stability} such non-FIFO allocation rule has been studied in the single-commodity case.
	\end{example}
	
	We conclude this section with the following standard result establishing well-posedness of the initial value problem associated to any MDFN. 
	
	\begin{lemma}\label{lemma:well-posedness}
		Consider a MDFN and define $\mc X$ be as in \eqref{X}. 
		Then, there exists a map 
		$$\varphi:\R_+\times\mc X\to\mc X\,,$$
		such that for every $\ov x$ in $\mc X$, $x(t)=\varphi(t,\ov x)$ is the unique solution of the MDFN with initial state $x(0)=\ov x$.
	\end{lemma}
	\ifshowproofs
	\begin{proof}
		See Appendix \ref{sec:proof-lemma-well-posedness}. 
	\end{proof}
	\fi
We conclude this section by recalling a few standard definitions that will be used in the next section. For a MDFN: 
	\begin{itemize} 
	\item for $\delta>0$ and $x^*$ in $\mc X$, 
	$$\mc B_{\delta}(x^*)=\{x\in\mc X:\,||x-x^*||_1<\delta\}\,;$$
	\item  an equilibrium point   is an array $x^*$ in $\mc X$ such that $\varphi(t,x^*)=x^*$ for every $t\ge0$; 
	\item an equilibrium point $x^*$ is stable if for every $\eps>0$ there exists some $\delta>0$ such that, for every $x$ in $\mc B_{\delta}(x^*)$, we have that $\varphi(t,x)\in\mc B_{\eps}(x^*)$ for every $t\ge0$;  
	\item the region of attraction of an equilibrium point $x^*$ is defined as  $\mc A(x^*)=\{x\in\mc X:\varphi(t,x)\stackrel{t\to+\infty}{\longrightarrow}x^*\}$;  
	\item an equilibrium point $x^*$ is asymptotically stable if it is stable and there exists $\delta>0$ such that $\mc B_{\delta}(x^*)\subseteq\mc A(x^*)$; 
	\item an equilibrium point $x^*$ is globally asymptotically stable if it is stable and $\mc A(x^*)=\mc X$. 
	\end{itemize}

	\section{Existence and asymptotic stability of free-flow equilibrium points}
	In this section we study the existence and asymptotic stability of  free-flow equilibrium points of MDFNs, using single-commodity studies in \cite{Coogan} as a foundation.

	Towards this goal, we start by introducing the notions of capacity region of a cell and of stability region of  an MTN. 
	\begin{definition}
		Consider a MTN with topology $\mc G=(\mc V,\mc E)$, set of commodities $\mc K$, supply functions $s_i$ and demand functions $d_i^{(k)}$ satisfying Assumption \ref{ass:demand+supply}, and  routing matrices $R^{(k)}$ satisfying Assumption \ref{ass:routing-matrix}, for every  $i$ in $\mc E$ and $k$ in $\mc K$. 
		Then: 
		\begin{enumerate}
			\item[(i)]	the \emph{capacity region} of a non on-ramp cell $i$ in $\mc E\setminus\mc R$  is 
			\begin{equation}\label{eq:cap_reg}
				\!\!\!\!\!\!\!\!\!\mathcal{C}_i = \left\{\zeta \in \mathbb{R}_+^\mathcal{K}: \sum\limits_{k \in \mathcal{K}}{\zeta_k} < s_i(\sum\limits_{k \in \mathcal{K}}{(d_{i}^{(k)})^{-1}(\zeta_k)})\right \}\,;
			\end{equation}
			\item[(ii)] the \emph{stability region} of the MTN is the set of inflow vectors $\lambda$ in $\mathbb{R}_+^{\mathcal{E}\times\mathcal{K}}$ such that 
			\begin{equation}\label{eq:stab_reg}
				\!\!\!\!\!\!\!\!\!\!\!\!\left(\left((I - (R^{(k)})^{\top})^{-1}\lambda^{(k)}\right)_i\right)_{k\in\mc K}\!\! \in \mathcal{C}_i\,,\qquad \forall i \in \mathcal{E}\setminus\mc R\,,
			\end{equation}
			and will be denoted by $$\Lambda=\left\{\lambda\in\mathbb{R}_+^{\mc E\times\mc K}:\,\eqref{eq:stab_reg}\right\}\,.$$ 
		\end{enumerate}
	\end{definition}
	A few comments are in order. First, consider the single-commodity case, i.e., the special case when $|\mc K|=1$. In this case, the capacity region of each non on-ramp cell $i$  in $\mc E\setminus\mc R$ is a semi-open interval $\mc C_i=[0,c_i)$ where the capacity value $$c_i=\sup\{d_i(\xi):\,\xi\ge0,\,  d_i(\xi)<s_i(\xi)\}\,,$$ is either the value where the cell's demand and supply curves meet, or it is equal $+\infty$, when the two curves do not meet. On the other hand, for every exogenous flow vector $\lambda$ in $\R_+^{\mc E}$, 
	\be\label{eq:transported}(I - R^{\top})^{-1}\lambda=\lambda+R^{\top}\lambda+(R^{\top})^2\lambda+(R^{\top})^3\lambda+\ldots\,,\ee
	is a nonnegative vector that accounts for both direct and indirect effects of $\lambda$, as transported by the routing matrix  $R$. 
	Hence, in the single-commodity case, \eqref{eq:stab_reg} requires that every entry $\left((I - R^{\top})^{-1}\lambda\right)_i$  of the vector \eqref{eq:transported} corresponding to a non on-ramp cell $i$ in $\mc E\setminus\mc R$ is strictly below the corresponding capacity value $c_i$. 
	
	In the multi-commodity case, i.e., when $|\mc K|>1$, \eqref{eq:cap_reg} defines the capacity region $\mc C_i$ of a non on-ramp cell $i$ in $\mc E\setminus\mc R$ as the set of non-negative vectors $\zeta$ in $\R^{\mc K}_+$ whose entry sum $\sum_k\zeta_k$ is strictly below the value of the cell's supply function $s_i(\cdot)$ computed in the sum across all commodities $k$ in $\mc K$ of the inverse demand functions $(d_i^{(k)})^{-1}(\zeta_k)$. For every exogenous inflow vector $\lambda^{(k)}$ in $\R_+^{\mc E}$ of a commodity $k$ in $\mc K$, each nonnegative vector $(I - (R^{(k)})^{\top})^{-1}\lambda^{(k)}$  again accounts for both direct and indirect effects of $\lambda^{(k)}$, as transported by commodity $k$'s routing matrix $R^{(k)}$. Condition \eqref{eq:stab_reg} requires that, for every non on-ramp cell $i$ in $\mc E\setminus\mc R$, the nonnegative vector $\zeta$ in $\R_+^{\mc K}$ with entries $\left((I - (R^{(k)})^{\top})^{-1}\lambda\right)_i$ for every commodity $k$ in $\mc K$, belongs to the capacity region $\mc C_i$. 
	
	Notice that, if  the demand functions $d_i^{(k)}$ of a non onramp cell $i$ in $\mc E\setminus\mc R$ are concave, their inverses $(d_i^{(k)})^{-1}$ are convex; if, moreover, the supply function $s_i$ is concave non-increasing, then the composition $\zeta\mapsto s_i(\sum_{k \in \mathcal{K}}{(d_{i}^{(k)})^{-1}(\zeta_k)})$ is concave, hence in this case the capacity region $\mc C_i$ of the cell $i$ is convex. Since \eqref{eq:stab_reg} dictates that a linear function of the exogenous flow $\lambda$ belongs to the capacity region $\mc C_i$ for every non onramp cell $i$ in $\mc E\setminus\mc R$, we get that a MTN whose demand and supply functions are all concave has a convex stability region $\Lambda$.
	
	\begin{contexample}
		Recall the demand and supply functions previously defined. Here we will set $\beta = \alpha = 1$ and $\gamma = 2$, which yields
		$
		d_i^{(k)}(\xi) = \xi 
		$
		and
		$
		s_i(\xi) = 2 - \xi
		$.
		Hence for every $\zeta_k$ in $\R_+$, this implies that
		$
		(d_i^{(k)})^{-1}(\zeta_k) = \zeta_k 
		$
		.
		For the sake of simplicity we shall assume to have only two commodities, namely \textit{a} and \textit{b}. Then, for any non-onramp cell $i$ in $\mc{E}\setminus\mc{R}$, we can compute the capacity region as
		$
		\mathcal{C}_i = \left\{\zeta \in \mathbb{R}_+^2: \zeta_a + \zeta_b < 1)\right \}\,.
		$
		The capacity region of a cell $i$ just found can be seen in Figure \ref{fig:cap_reg}.
		\begin{figure}
			\centering
			\includegraphics[width=0.3\linewidth]{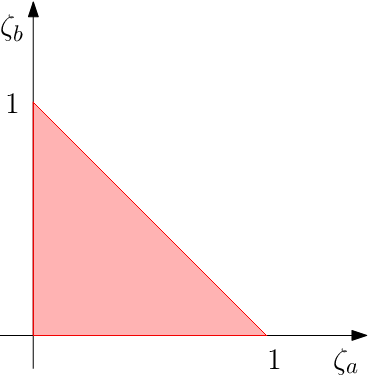}
			\caption{Capacity region with $|\mc K|=2$, $d_i^{(1)}\!(\xi) =d_i^{(2)}\!(\xi) = \xi$, $s_i(\xi) = 2 - \xi$.}
			\label{fig:cap_reg}
		\end{figure}
	\end{contexample}

	
	We are now ready to state our first result providing necessary and sufficient conditions for the existence and uniqueness of a free-flow equilibrium in a MDFN. 
	\begin{proposition}
		Consider a MTN satisfying Assumptions \ref{ass:demand+supply} and \ref{ass:routing-matrix}, and let $\Lambda$ be its stability region. 
		Then, every MDFN with exogenous inflow array $\lambda$ admits a free-flow equilibrium point if and only if $\lambda\in\Lambda$. Moreover, in this case, the free-flow equilibrium point is unique. 
	\end{proposition}
	\ifshowproofs
	\begin{proof}
		It follows from \eqref{compact} that $x^*$ in $\mc X$  is an equilibrium point for the MDFN if and only if 
		\begin{equation}\label{eq:inverse}
			z^{(k)}(x^*) = (\mathcal{I}-(R^{(k)})^{\top})^{-1}\lambda^{(k)}\,,\qquad \forall k \in \mathcal{K}\,.
		\end{equation}

		Now, if $x^*$ in $\mc F$ is a free-flow equilibrium point, then  it follows from \eqref{eq:free_dyn} and \eqref{eq:zik} that 
		\be\label{zx}z_i^{(k)}(x^*)=d_{i}^{(k)}((x^*)_{i}^{(k)})\,,\ee
		for every cell $i$ in $\mathcal{E}$ and commodity $k$ in $\mc K$. 
		Equations \eqref{zx}  and \eqref{free-flow} imply  that 
		$$\ba{rcl}\sum\limits_{k \in \mathcal{K}}z_{i}^{(k)}(x^*)
		&=&\sum\limits_{k \in \mathcal{K}}{\sum\limits_{j \in \mathcal{E}}{R_{ji}^{(k)}}d_{j}^{(k)}((x^*)_{j}^{(k)})}  \\
		&<&  s_i(\sum\limits_{k \in \mathcal{K}}{(x^*)_{i}^{(k)}})\\
		&=& s_i(\sum\limits_{k \in \mathcal{K}}(d_{i}^{(k)})^{-1}(z_i^{(k)}(x^*)))\,,\ea$$
		i.e., $(z_{i}^{(k)}(x^*))_{k\in\mc K}\in\mc C_i$ for every non on-ramp cell $i$ in $\mc E\setminus\mc R$. Along with \eqref{eq:inverse}, this implies that $\lambda\in\Lambda$. 


		
		Conversely, define the arrays $\zeta$ and $x^*$ in $\mc X$ with entries
		\be\label{zeta}\zeta_i^{(k)}= \left((\mathcal{I}-(R^{(k)})^{\top})^{-1}\lambda^{(k)}\right)_i\,,\ee
		\be\label{x*}(x^*)^{(k)}_i=(d_i^{(k)})^{-1}(\zeta^{(k)}_i)\,,\ee
		respectively, for every $ i$ in $\mc E$  and $ k$ in $\mc K$.
		If  $\lambda\in\Lambda$, then $(\zeta_i^{(k)})_{k\in\mc K} \in \mathcal{C}_i$ for every non on-ramp cell $i$ in $\mc E\setminus\mc R$, so that
		$$\sum\limits_{k \in \mathcal{K}}d_{i}^{(k)}((x^*)_{i}^{(k)})=\sum\limits_{k \in \mathcal{K}}{\zeta^{(k)}_i} < s_i(\sum\limits_{k \in \mathcal{K}}((x^*)_{i}^{(k)})\,,$$
		i.e., $x^*\in\mc F$ is in free-flow. 
		It then follows from \eqref{eq:free_dyn} that 
		$$z_{i}^{(k)}(x^*) = d_{i}^{(k)}((x^*)_{i}^{(k)})=\zeta^{(k)}_i=\left((\mathcal{I}-(R^{(k)})^{\top})^{-1}\lambda^{(k)}\right)_i\,,$$
		for every cell $i$ in $\mc E$ and commodity $k$ in $\mc K$, so that \eqref{eq:inverse} is satisfied. Hence $x^*$ is a free-flow equilibrium. Finally, notice that the entries of every free-flow equilibrium $x^*$ must necessarily be the ones determined by \eqref{zeta}--\eqref{x*}, so that the free-flow equilibrium $x^*$ is unique, whenever it exists. 
		%
	\end{proof}
	\fi

	We shall now introduce the concept of $l_1$-nonexpansivness useful for the stability analisys.

	\begin{definition}
		The MDFN is $l_1$-nonexpansive in a region $\mc{R}\subseteq\mc X$ if 
		for every $x$ and $y$ in $\mc R$ and for every $T\ge0$ such that 
		\begin{equation}\varphi(s,x),\varphi(s,y) \in \mathcal{R}\,, \qquad \forall \; 0 \leq s \leq T\,,
		\end{equation}
		we have that 
		\begin{equation}
			||\varphi(t,y) - \varphi(t,x)||_1 \leq ||y - x||_1
		\end{equation}
	\end{definition}
	\begin{lemma}\label{lemma:non-expansive}
		Every MDFN is $l_1$-nonexpansive in every convex subset $\mc D\subseteq\mc F$ of its free-flow region.
	\end{lemma}
	\ifshowproofs
	\begin{proof} See Appendix \ref{sec:proof-lemma-non-expansive}.\end{proof}
	\fi

	\begin{lemma}\label{lemma:1}
		Consider a MDFN. 
		For a free-flow equilibrium point $x^*$  in $\mathcal{F}$, let 
		\begin{equation*}
			\delta^* = \sup \{\delta>0: \text{MDFN is }l_1\text{-nonexpansive on }\mathcal{B}_{\delta}(x^*)\}
		\end{equation*}
		If $x^*$ is asymptotically stable, then $\mathcal{B}_{\delta^*}(x^*)\subseteq\mathcal{A}(x^*)$.
	\end{lemma}
	
	\ifshowproofs
	\begin{proof}See Appendix \ref{sec:proof-lemma-1}. \end{proof}
	\fi
	
	We are now ready to state and prove our main result. 
	
	\begin{theorem}\label{theorem:main}
		Consider a MTN satisfying Assumptions \ref{ass:demand+supply} and \ref{ass:routing-matrix} and let $\lambda$ in $\Lambda$ be an exogenous inflow array within its stability region. 
		Consider a MDFN and let $x^*$ in $\mc F$ be its unique free-flow equilibrium. 
		Define 
		\begin{equation}
			\bar{\delta} = \underset{x \in \mathcal{X} \setminus \mathcal{F}}{\inf} ||x - x^*||_1 > 0\,.
		\end{equation}
		Then, $x^*$ is asymptotically stable and 
		$$\mathcal{B}_{\bar{\delta}}(x^*) \subseteq \mathcal{B}_{\delta^*}(x^*) \subseteq \mathcal{A}(x^*)\,.$$
	\end{theorem}
	\ifshowproofs
	\begin{proof}
		It is sufficient to prove that $x^*$ is asymptotically stable, as the rest of the claim follows from Lemmas \ref{lemma:non-expansive} and \ref{lemma:1}. 
We first asymptotic stability of $x^*$. For every commodity $k$ in $\mc K$, let
$$				f^{(k)}(x^{(k)})  = \lambda^{(k)} - (I - (R^{(k)})^{\top})z^{(k)}(x)\,.$$
Then, 		\begin{equation}
				A^{(k)}=\nabla_{x^{(k)}} f^{(k)}((x^*)^{(k)})  =  (I - (R^{(k)})^{\top})D^
			{(k)}\,,
		\end{equation}
		where $D^{(k)}$ in $\R^{\mc E\times\mc E}$ is a diagonal matrix with positive diagonal elements $(d^{(k)}_i)'((x^*)^{(k)}_i)>0$ for every $i$ in $\mc E$. 
Notice that the matrix $A^{(k)}$ is compartmental, i.e., $\sum_jA^{(k)}_{ij}\le0$ for every $i$ in $\mc E$. By the Perron-Frobenius Theory, there exists a real eigenvalue $\lambda^{(k)}$ in $\R$ and a nonnegative eigenvector $y$ in $\mathbb{R}^\mathcal{E}_+$ such that $(A^{(k)})y = \lambda^{(k)}y$ and $\Re(\lambda)\le\lambda^{(k)}$ for every other eigenvalue $\lambda$ of $A^{(k)}$. Then, define 		\begin{equation}\label{eq:schur}
				\mathcal{J} = \{ i = 1, \dots ,n : y_i > 0\}, \qquad y_* = \underset{j \in \mathcal{J}}{\min}y_j > 0
				\end{equation}
 and let $d_* = \min \{D_{ii}^{(k)} : i \in \mathcal{J} \}$. Since the matrix $R^{(k)}$ is sub-stochastic and out-connected, we have that
$$\gamma=\min_{i\in\mc J}\sum\limits_{j \in \mathcal{J}}R_{ij}^{(k)}<1\,,$$
so that
		\begin{equation}
			\underset{i \in \mathcal{J}}{\min}\sum\limits_{j \in \mathcal{J}}{A_{ij}^{(k)}} = \underset{i \in \mathcal{J}}{\min}{D_{ii}^{(k)}(1 - \sum\limits_{j \in \mathcal{J}}R_{ij}^{(k)})} \leq -d_*\gamma < 0\,.
		\end{equation}
		Then, we have
		\begin{equation}
			\begin{array}{rcl}
				\lambda^{(k)}\sum\limits_{j \in \mathcal{J}}{y_j} & = & \sum\limits_{i\in\mc E}{\sum\limits_{j \in \mathcal{J}}{A_{ij}y_i}} \\ & = &\sum\limits_{i \in \mathcal{J}}{y_i}\sum\limits_{j \in \mathcal{J}}{A_{ij}} \leq  -y_*d_*\gamma < 0
			\end{array}
		\end{equation}
		thus proving that $\lambda^{(k)}< 0$. Hence, the matrix $A^{(k)}$ is Hurwitz-stable for every commodity $k$ in $\mc K$, yielding that the free-flow equilibrium point $x^*$ is indeed locally asymptotically stable. 
	\end{proof}
\fi
We now show that, as a special case, Theorem \ref{theorem:main} recovers a result first proven in \cite{Stability}. 
	\begin{corollary}
Consider a MTN with $|\mc K|=1$, satisfying Assumptions \ref{ass:demand+supply} and \ref{ass:routing-matrix} and let $\lambda$ in $\Lambda$ be an exogenous inflow array within its stability region. 
		Consider the MDFN with the non-FIFO allocation rule \eqref{eq:10}. 		
		Then, the free-flow equilibrium point $x^*$ is globally asymptotically stable.
	\end{corollary}
	\ifshowproofs
	\begin{proof}
		In a single-commodity scenario, define the system's dynamic as
		\begin{equation}
			\dot{x} = \lambda - (\mathcal{I} - R^{\top})z(x) = \phi(x)
		\end{equation}
		Notice that the system is monotone since its gradient
		\begin{equation}
			\frac{\partial \phi_i}{\partial x_j}(x) \geq 0 \quad \forall i \neq j
		\end{equation}
		is Metzler. Moreover, notice that
		\begin{equation}
			\sum\limits_{j \in \mathcal{E}}{\frac{\partial \phi_i}{\partial x_j}(x)} \leq 0 \quad \forall i \in \mathcal{E}
		\end{equation}
		As we have already proven, this means that the system is $l_1$-nonexpansive on the whole state space $\mathcal{X}$.
		Thus, we have $\delta^* = +\infty$ and $\mathcal{B}_{\delta^*}(x^*) = \mathcal{A}(x^*) = \mathcal{X}$.
	\end{proof}
	\fi
	The following result concerns another special case when Theorem \ref{theorem:main} ensures global asymptotic stability of the free-flow equilibrium point. 
	
	\begin{corollary}
		Given $\lambda \in \Lambda$, let $x^* \in \mathcal{F}$ be such that 
		\begin{equation}
			(x_{i}^{(k)})^* = (d_{i}^{(k)})^{-1}([(\mathcal{I}-(R^{(k)})^{\top})^{-1}\lambda^{(k)}]_i) \quad \forall i\in \mathcal{E}, k \in \mathcal{K}
		\end{equation}
		If for every $t \geq 0$ holds that
		\begin{equation}
			s_i(\sum\limits_{k \in \mathcal{K}}{x_{i}^{(k)}}) > \sum\limits_{k \in \mathcal{K}}{d_{i}^{(k)}(x_{i}^{(k)})} \quad \forall i \in \mathcal{E}
		\end{equation}
		then the equilibrium point $x^*$ is GAS.
	\end{corollary}
	\ifshowproofs
	\begin{proof}
		Notice that if the supply function is always larger than the sum of the demand functions of the different commodities, then the MDFN is always in free-flow. Thus, we have
		\begin{equation}
			\delta^*  = \bar{\delta} = \underset{x \in \mathcal{X} \setminus \mathcal{F}}{\inf} ||x-x^*||_1 = +\infty
		\end{equation}
		yielding that
		\begin{equation}
			\mathcal{B}_{\bar{\delta}}(x^*) = \mathcal{B}_{\delta^*}(x^*) = \mathcal{A}(x^*) = \mathcal{X}\,,
		\end{equation}
		so that the equilibrium point $x^*$ is globally asymptotically stable.
	\end{proof}
	\fi
	\section{An Analytical Example}
	Even though it is possible to characterize free-flow equilibrium points analytically, it does not appear clear whether some non free-flow equilibrium points may arise and, if they do so, under what circumstances.
	In order to address this problem, we introduce a relevant example, while keeping the dimensionality of the MDFN as low as possible.
	\par
\begin{figure}
		\centering
		\includegraphics[width=0.5\linewidth]{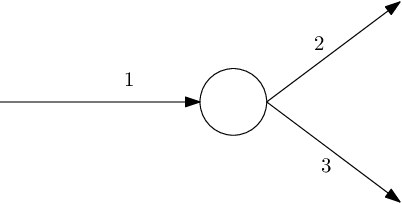}
		\caption{Diverge junction network}
		\label{fig:ex_cong}
	\end{figure}

	Consider the network $\mathcal{G} = (\mathcal{V},\mathcal{E})$ in Figure \ref{fig:ex_cong} with $\{1\}$ in $\mathcal{R}$ and $\{2,3\}$ in $\mathcal{S}$. The flow is split into two commodities, namely \textit{a} and \textit{b}, that interact with each other in a non-FIFO manner. Assume demand and supply functions as in Example \ref{ex:1}, where for each demand  $\xi = x_i^{(k)}$ and for each supply $\xi = \sum_k{x_i^{(k)}}$. Commodity \textit{a} is equally split between cells ${2}$ and ${3}$, while commodity \textit{b} is sent entirely towards cell $2$. Thus, the routing matrices, $R^{(a)}$ and $R^{(b)}$, can be defined such that $R^{(a)}_{12} = R^{(a)}_{13} = 0.5$, $R^{(b)}_{12} = 1$ and equal to 0 otherwise.
	\par
	Once our framework has been set, we can delve into the study of the stability of the MDFN. Firstly, it is crucial to understand under what conditions an equilibrium does exist, i.e., the maximum amount of exogenous inflow that can be handled by the network.
	To this end, the fact that at equilibrium the difference between the total incoming and out-going flow must be zero implies the existence of a region $\Lambda_B$, defined as 	
	$
		\Lambda_{B} = \{\lambda \in \mathbb{R}^2 : \lambda_{1}^{(a)}+ \lambda_{1}^{(b)} \leq 2 \}
	$.
	Recall that, in $\eqref{eq:stab_reg}$, we defined the region for which a free-flow equilibrium point does exist. In our example, this stability region is
	$
		\Lambda_{FF} := \Lambda = \{\lambda \in \mathbb{R}^2 : \lambda_{1}^{(a)} + 2\lambda_{1}^{(b)}< 2\}
	$.
	Notice that $\Lambda_{FF} \subseteq \Lambda_B$, which follows from the fact that $\Lambda_B$ is the biggest region that makes an equilibrium admissible. Moreover, this implies the existence of a third region $\Lambda_{B \setminus FF}$ such that $\Lambda_{FF} \cup \Lambda_{B \setminus FF} = \Lambda_B$.
	
	Only four scenarios may arise: (i) free-flow, (ii) congestion on cell $2$, (iii) congestion on cell $3$, and (iv) congestion on both cells $2$ and $3$.
	Our analysis will then proceed by analytically finding the equilibrium points in each of those. 
	
	\begin{enumerate}
		\item[(i)] \textbf{free-flow:} as proven, a free-flow equilibrium point exists only if $\lambda$ is in $\Lambda_{FF}$ and can be computed as
		\begin{equation*}
			\begin{array}{rcl}
				(x^{(a)})^* & = & (\mathcal{I} - (R^{(a)})^\top)^{-1}\lambda^{(a)} \\ (x^{(b)})^* & = & (\mathcal{I} - (R^{(b)})^\top)^{-1}\lambda^{(b)}
			\end{array}
		\end{equation*}
		Thus, we find
		\begin{equation*}
			\begin{array}{rcl}
				(x^{(a)})^* & = & (\lambda_{1}^{(a)},\; \frac{1}{2}\lambda_{1}^{(a)}, \; \frac{1}{2}\lambda_{1}^{(a)}) \\
				(x^{(b)})^* & = & (\lambda_{1}^{(b)}, \; \lambda_{1}^{(b)})
				
			\end{array}
		\end{equation*}
		\item[(ii)] \textbf{Congestion only on cell 2:} in the non-FIFO model this implies that
		\begin{equation}\label{eq:2_cong}
			0 < \frac{2 - x_{2}^{(a)} - x_{2}^{(b)}}{\frac{1}{2}x_{1}^{(a)} + x_{1}^{(b)}} < 1
		\end{equation}
		Since $x_{2}^{(a)} + x_{2}^{(b)} = 1$, because cell 2 is a congested off-ramp, we solve the equations of equilibrium, yielding
		
		\begin{equation*}
			\begin{array}{rcl}
				(x^{(a)})^* & = & (2(\lambda_{1}^{(a)}, \; 1-\lambda_{1}^{(b)}, \; \lambda_{1}^{(a)} + \lambda_{1}^{(b)} - 1) \\
				(x^{(b)})^* & = & (\frac{\lambda_{1}^{(b)}(\lambda_{1}^{(a)} + \lambda_{1}^{(b)} - 1)}{1 - \lambda_{1}^{(b)}}, \; \lambda_{1}^{(b)})
			\end{array}
		\end{equation*}
		Hence, the equilibrium point is unique and by substituting it in \eqref{eq:2_cong} yields the inequality
		\begin{equation*}
			\lambda_{1}^{(a)} + 2\lambda_{1}^{(b)} > 2 \rightarrow \lambda \in \Lambda_{B \setminus FF}
		\end{equation*}
		So, a congestion happening on cell 2 allows for an equilibrium $x^*$ only if $\lambda$ is chosen outside of the stability region.
		\item[(iii)] \textbf{Congestion only on cell 3:} in this particular setting this case never occurs. Indeed, since cell 3 always receives less flow, if it is congested so is cell 2;
		\item[(iv)] \textbf{Congestion on both cells 2 and 3:} this case can be interpreted as
		\ben
			\ba{rcl}
			0  < & \frac{2-x_{2}^{(a)}-x_{2}^{(b)}}{\frac{1}{2}x_{1}^{(a)}+x_{1}^{(b)}} & < 1 \\
			0  < & \frac{2-x_{3}^{(a)}}{\frac{1}{2}x_{1}^{(a)}} & < 1
			\ea
		\een
		Again, it is possible to solve the equilibrium's equation, yielding
		\begin{equation*}
			\begin{array}{rcl}
				(x^{(a)})^* & = & (c\displaystyle\frac{\lambda_{1}^{(a)} - 1}{\frac{1}{2}\lambda_{1}^{(b)}}, \; \lambda_{1}^{(a)} - 1, \; 1) \\
				(x^{(b)})^* & = & (c, \; \lambda_{1}^{(b)})
			\end{array}
		\end{equation*}
		where $c$ is in $\mathbb{R}_+$. So, there are infinitely many equilibrium points, that have different values of $x_{1}^{(a)}$ and $x_{2}^{(b)}$. Moreover, if we substitute an equilibrium point into the condition $x_{2}^{(a)} + x_{2}^{(b)} = 1$, we obtain
		$
			\lambda_{1}^{(a)} + \lambda_{1}^{(b)} = 2$ so that $ \lambda \in \partial \Lambda_B
		$.
		To summarize, a continuum of equilibrium points can arise if $\lambda$ is chosen on the boundary of $\Lambda_B$.
	\end{enumerate}
	\begin{figure}
		\centering
		\includegraphics[width = 0.3\linewidth]{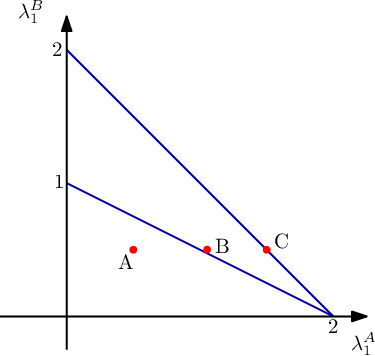}
		\caption{Stability and bounded regions with points $A=(0.5,0.5)$, $B=(1.2,0.5)$ and $C=(1.5,0.5)$}
		\label{fig:stab}
	\end{figure}
	To validate these results, we provide simulation of the system by choosing three different combinations of exogenous inflows as shown in Figure \ref{fig:stab}. Moreover, for each combination two different initial conditions $x(0)$ and $\tilde{x}(0)$ are considered.
	\par
	Given the exogenous inflows $\lambda$, it is possible to compute the equilibrium points for each case separately.
	\ben
		\begin{array}{rcl}
			x^*_{FF} & = & (0.5, 0.5, 0.25, 0.5, 0.25)\\
			x^*_{B \setminus FF} & = & (1.4, 0.7, 0.5, 0.5, 0.7)\\
			x^*_{\partial B} & = & (2c, c, 0.5, 0.5, 1)
		\end{array}
	\een
	\par
\begin{figure}
		\begin{minipage}{0.16\textwidth}
			\centering
			\includegraphics[width=1.\linewidth]{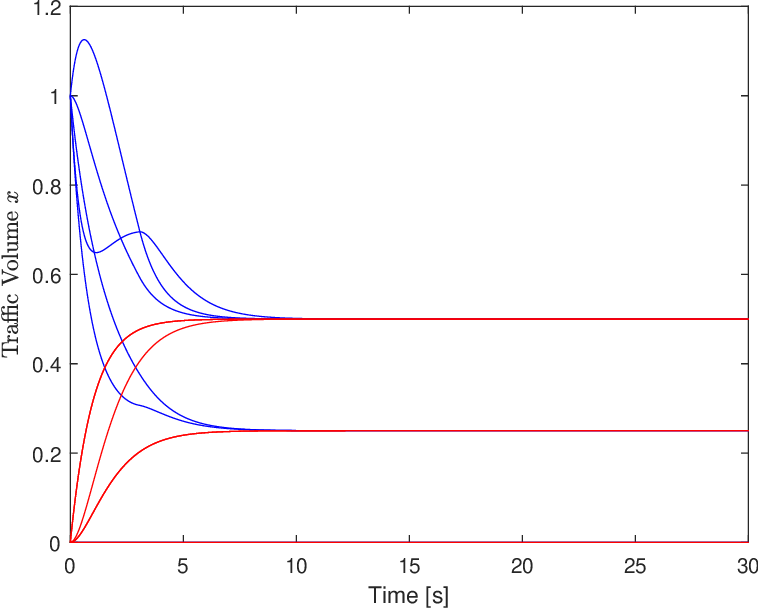}
		\end{minipage}\hfill
		\begin{minipage}{0.16\textwidth}
			\centering
			\includegraphics[width=1.\linewidth]{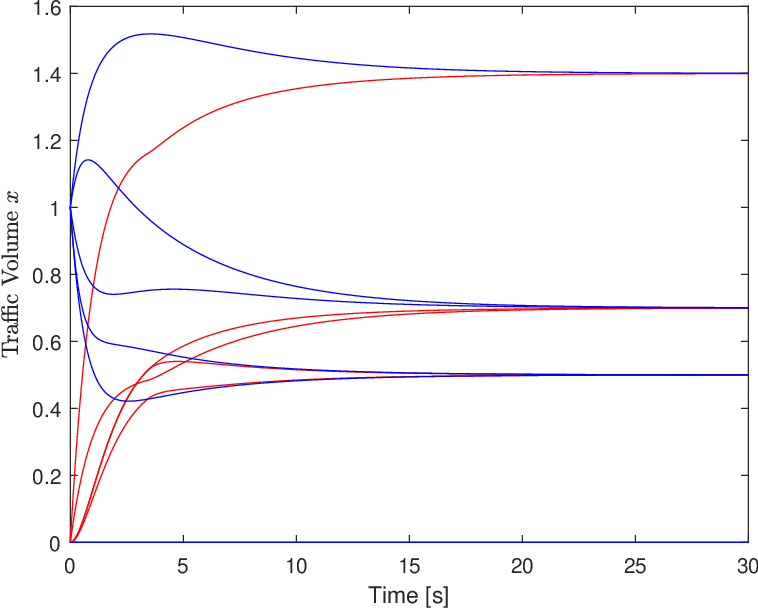}
		\end{minipage}\hfill
		\begin{minipage}{0.16\textwidth}
			\centering
			\includegraphics[width=1.\linewidth]{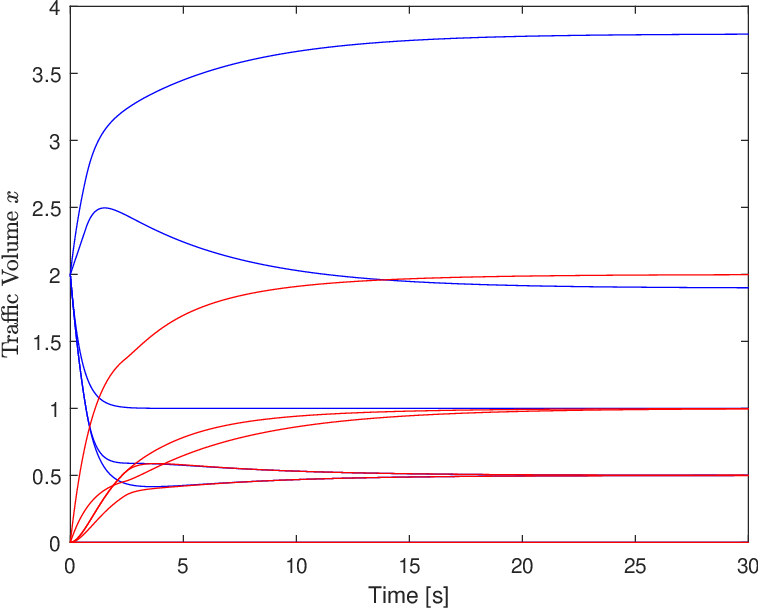}
		\end{minipage}
		\caption{Two initial conditions ($x(0)$ in red and $\tilde{x}(0)$ in blue) under different exogenous inflows: $\lambda \in \Lambda_{FF}$ (left), $\lambda \in \Lambda_{B \setminus FF}$ (middle) and $\lambda\in \Lambda_B$ (right).}
		\label{fig:results}
	\end{figure}

	The results, shown in Figure \ref{fig:results}, show that if $\lambda$ is chosen inside or outside the stability region, the trajectories corresponding to the two initial conditions converge to the same equilibrium point. On the contrary, when $\lambda$ is chosen on the boundary of $\Lambda_{B}$, different initial condition may converge to different equilibrium points.
	
	\section{Conclusions}
	This paper studies a class of dynamical systems based on dynamical flow networks led by mass conservation laws, that act similarly in the free-flow region. This allows us to analyze its stability properties and discuss the presence of equilibrium points. Moreover, through an insightful example it is shown the existence of equilibrium points not in free-flow. Future researches could focus on generalizing the presence of congested equilibrium points to any network and on studying the system's global stability.

	\bibliographystyle{ieeetr}
	\bibliography{references}
	
	\ifshowproofs
	\appendices
	\section{Proof of Lemma \ref{lemma:well-posedness}}\label{sec:proof-lemma-well-posedness}
	%
	
	Notice that this is always true based on the assumption that if the initial conditions $x(0)$ have all nonnegative entries, so does $x(t)$ for all $t \geq 0$. A sufficient condition for this is the assumption made on $d_{i}^{(k)}(x_{i}^{(k)})$, which is $d_{i}^{(k)}(x_{i}^{(k)}) = 0$ whenever $x_{i}^{(k)} = 0$.
	\qed 
	
	\section{Proof of Lemma \ref{lemma:non-expansive}}\label{sec:proof-lemma-non-expansive}
	Define the $l_1$ matrix measure of the gradient $\nabla f(x)$ as
$$
		\begin{array}{rcl}			
			\mu(\nabla f(x)) & = &\underset{j \in \mathcal{E}}{\max} \sum\limits_{i \in \mathcal{E}}{\frac{\partial f_i}{\partial x_j^{(k)}}(x)}\\ & = & \underset{j \in \mathcal{E}}{\max} \sum\limits_{i \in \mathcal{E}}{R_{ji}^{(k)}(z_j^{(k)})'(x) - (z_j^{(k)})'(x)}\,. 
		\end{array}
$$
Note that $\mu(\nabla f(x))\le0$. 
	We shall prove that this implies that the system is $l_1$-nonexpansive. Let $\mc D\subseteq\mc F$ be nonempty convex. Fix $x,y$ in $\mc D$ and let
$\psi(t,\lambda) = \varphi(t,(1-\lambda)x + \lambda y)$,  $z_\lambda(t) = \frac{\partial \psi}{\partial \lambda}(t,\lambda)$, and  $A(t) = \nabla_x f(\psi(t,\lambda))$.	
Observe that
$$
		\begin{array}{rcl}
			\dot{z}_\lambda(t) & = & \frac{\partial}{\partial t}\frac{\partial \psi}{\partial \lambda}(t,\lambda) = \frac{\partial}{\partial \lambda}\frac{\partial \psi}{\partial t}(t,\lambda) = \frac{\partial}{\partial \lambda} f(\psi(t,\lambda)) \\
			& = & \nabla_x f(\psi(t,\lambda)) \frac{\partial \psi}{\partial \lambda}(t,\lambda) = A(t)z_\lambda(t)
		\end{array}
$$
	So, it is possible to use Coppel Inequality to find that
$$
		||z_\lambda(t)||_1 \leq \exp(\int_0^t(\mu(A(s)ds)))||z_\lambda(0)||_1 \leq ||y-x||_1\,.
$$
	Then from the Fundamental Theorem of Calculus, we get
$$
		\varphi(t,y) - \varphi(t,x) = \int_0^1 \frac{\partial \psi}{\partial \lambda}(t,\lambda)d\lambda = \int_0^1 z_\lambda(t) d\lambda\,,
$$
	which leads to
$$
		\begin{array}{rlc}
			||\varphi(t,y) - \varphi(t,x)||_1 & = & ||\int_0^1 \frac{\partial \psi}{\partial \lambda}(t,\lambda)d\lambda||_1 \\ & \leq & \int_0^1 ||z_\lambda(t)||_1d\lambda \leq ||y-x||_1\,. 
		\end{array}
$$
	So the system is $l_1$-nonexpansive in the free-flow region. \qed
	
	\section{Proof of Lemma \ref{lemma:1}} \label{sec:proof-lemma-1}
	This proof will follow the reasoning used in \cite[Lemma~6]{Stability}. Sufficiency is obvious, so we shall evaluate only necessity.
	\par
	Given $x^*$ a LAS equilibrium point, then there exists a $\mathcal{KL}$ function $\beta(\cdot,\cdot)$ such that
	$||\varphi(t,x) - x^*||_1 \leq \beta(x-x^*,t)$ for every $x$ inside a sufficiently small closed ball $\mathcal{B}_\varepsilon(x^*)$ of radius $\varepsilon$ in the $l_1$-topology. This means that for every $x^\circ \in \overline{\mathcal{B}_\varepsilon(x^*)}$ then $\varphi(t,x^{\circ}) \rightarrow x^*$. So, we shall now prove that for every $x \in \mathcal{B}_{\delta^*}(x^*)$ there exists a T $\geq 0$ such that $\varphi(T,x) \in \mathcal{B}_\varepsilon(x^*)$. Let $\tilde{x} = x^* + \frac{\varepsilon}{||x - x^*||_1}(x-x^*)$, from which follows that $||\tilde{x}-x^*||_1 = \varepsilon$ and $||x - x^*||_1 = ||x-\tilde{x}||_1 + ||\tilde{x} - x^*||_1 = ||x - \tilde{x}||_1 + \varepsilon$. The fact that the MDFN is $l_1$-nonexpansive on $\mathcal{B}_{\delta^*}(x^*)$ implies that it is also positively invariant, so it must hold that $||\varphi(t,x) - \varphi(t,\tilde{x})||_1 \leq ||x-\tilde{x}||_1$. The triangle inequality implies that 
	\begin{equation*}
		\begin{array}{rcl}
			\!\!||\varphi(t,x) - x^*||_1 \!\!\!&\!\!\leq\!\!&\!\! ||\varphi(t,x) - \varphi(t,\tilde{x})||_1 + ||\varphi(t,\tilde{x}) - x^*||_1 \\
			&\!\! =\!\! &\!\! ||x - \tilde{x}||_1 + ||\varphi(t,\tilde{x}) - x^*||_1 \\
			&\!\! =\!\! &\!\! ||x - x^*||_1 - \varepsilon + ||\varphi(t,\tilde{x}) - x^*||_1
		\end{array}
	\end{equation*}
	Since $\beta(\cdot,\cdot) $ is a $\mathcal{KL}$ function, there exists $T_{\varepsilon} \geq 0$ such that $\beta(x-y,t) \leq \frac{\varepsilon}{2}$ for every $y$ in $\overline{\mathcal{B}_{\varepsilon}(x^*)}$ and for $t \geq T_{\varepsilon}$. So,
	\begin{equation*}
		\begin{array}{rcl}
			||\varphi(t,x) - x^*||_1 &\leq& ||x - x^*||_1 - \varepsilon + ||\varphi(t,\tilde{x}) - x^*||_1 \\
			& \leq &|x - x^*||_1 - {\varepsilon}/{2}
		\end{array}
	\end{equation*}
	for every $t \geq T_{\varepsilon}$. If $\varphi(T_{\varepsilon}, x) \in \overline{\mathcal{B}_{\varepsilon}(x^*)}$, then the proof is complete. If this is not the case, the same argument can be iterated so that $\varphi(T,x) \in \mathcal{B}_{\varepsilon}(x^*)$.\qed
	\fi
\end{document}